# Fractal trajectories of the dynamical system

*Marek Berezowski*
*Silesian University of Technology*
*Institute of Mathematics*
*Gliwice, Poland*
E-mail address: marek.berezowski@polsl.pl

**Abstract**
The scope of the paper is a theoretical analysis of the dynamical system, the model of which was reduced to Weierstrasse function. A fractal structure of the trajectory was proved and the entropy of the system information designated.

## 1. Introduction

The dynamics of a mechanical device consisting of *n* wheels of different radii and rotating at different speed was tested. On the rim of a given wheel rotating with specific speed another smaller wheel is mounted, rotating at a higher speed. (Fig.1).

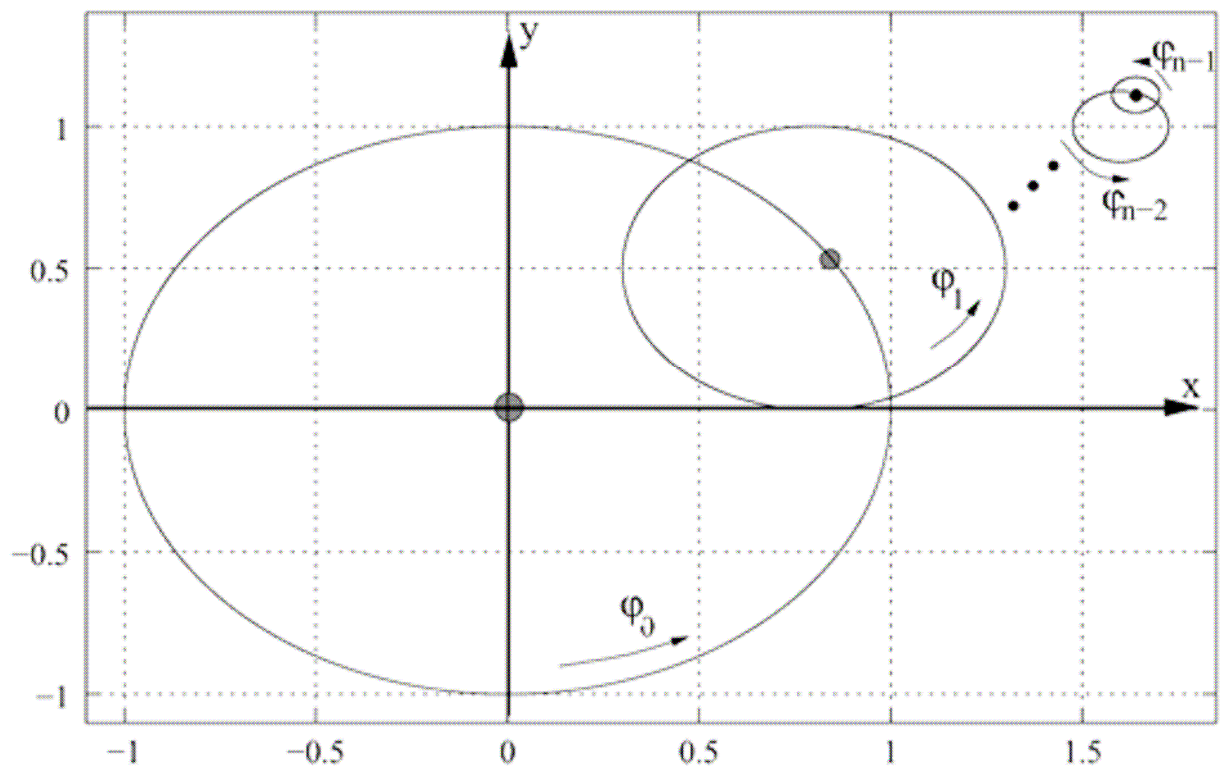

**Fig. 1.** Conceptual diagram of the system of wheels.

Accordingly, the particular wheels rotate not only at their own speed but are also driven by the motion of all the wheels that are bigger. On the rim of the last (the smallest) wheel a scriber is mounted, leaving traces on the plane. It turns out that even with a relatively small number of the wheels, despite a non-chaotic motion of the scriber, the trace may be very complex and may create a fractal figure [9, 10]. It is an interesting fact that the speed of the scriber is constantly changing, both in terms of its value and direction, although the angular speeds of the wheels are constant. In an extreme case, when the number of the wheels approaches infinity, the speed of the scriber changes at each moment. Such behavior is similar to Brownian motion.

## 2. Model of the system

As already mentioned in the Introduction, let us consider a mechanical device consisting of *n* different wheels. On the rim of the *k-th* wheel axis *(k+1)-th* is mounted (Fig. 1). The wheels rotate at specific speeds in the counter-clockwise direction. Let us assume that the rations between the radii and the angles of rotation of the wheels are constant and equal to:

$$q = \frac{r_k}{r_{k+1}} = \frac{\varphi_{k+1}}{\varphi_k} > 1. \tag{1}$$

The scriber mounted on the rim of the last wheel leaves traces on the surface. The position of the scriber is determined by its coordinates:

$$x' = \sum_{k=0}^{n-1} r_k \cos(\varphi_k) = r_0 \sum_{k=0}^{n-1} \frac{1}{q^k} \cos(q^k \varphi_0) \tag{2}$$

$$y' = \sum_{k=0}^{n-1} r_k \sin(\varphi_k) = r_0 \sum_{k=0}^{n-1} \frac{1}{q^k} \sin(q^k \varphi_0). \tag{3}$$

In terms of a complex notation, the above equations may be expressed as:

$$R' = \sum_{k=0}^{n-1} r_k e^{i\varphi_k} = r_0 \sum_{k=0}^{n-1} \frac{1}{q^k} e^{iq^k \varphi_0} \tag{4}$$

where $R'$ is the complex radius-vector.

On the grounds of the above relationships it is easy to indicate that for *n=2* and *q=2* the scriber marks the line the length of which is $L = 8r_0$. For $r_0 = 1/2$ the line demarcates the area of identical size and shape as the biggest Mandelbrot fractal, often referred to in literature as "*the heat curve*" [2].

By introducing the standardized coordinates and the reduced notation of the angular variable:

$$R = \frac{R'}{r_0}; \quad t = \varphi_0 \tag{5}$$

equation (4) is transformed into:

$$R = \sum_{k=0}^{n-1} \frac{1}{q^k} e^{iq^k t}. \tag{6}$$

This dependence is analogical to Weierstrasse function [5-6], the only difference being that the derivative of Weierstrasse function is infinite; whereas, in the discussed model, the derivatives are finite. As *R* changes in a continuous mode, the trajectory of the motion of the scriber is also continuous. However, the rate of the change of the radius is discontinuous for $n \to \infty$:

$$\frac{dR}{dt} = i \lim_{n \to \infty} \sum_{k=0}^{n-1} e^{iq^k t}. \tag{7}$$

For $n = \infty$ the discontinuity occurs for each value of angle $t$. To prove this, let us consider the following boundary:

$$\varepsilon = i \lim_{\Delta t \to 0} \sum_{k=0}^{\infty} \left( e^{i(q^k t + q^k \Delta t)} - e^{iq^k t} \right). \tag{8}$$

As index $k$ is directly equal to infinity, and the increase of angle $\Delta t$ only approaches zero, assuming that $q>1$ – quotient $q^k \Delta t$ is equal to infinity. This means that for each value of $t$ boundary $\varepsilon$ is different from zero, which proves that the derivative presented above is discontinuous for each value of $t$. In practice, when $n < \infty$ this phenomenon is important to $\Delta t > q^{1-n}$. Thus, for example, for $n=20$ and $q=2.5$, $\Delta t > 0.00000003$. This means that the loss of the discontinuity of the derivative is observable in the course of the measurements of the motion of the biggest wheel only at intervals shorter than 30 nano-degrees. Accordingly, in practice the system is subject of incessant rapid changes of the speed of the scriber motion, both in terms of its value and direction. It may be stated that the scriber is continuously driven, creating, in consequence, very complex fractal figures on the surface.

## 3. Calculations and analysis of the results

To set an example, two cases were considered: $q=2.5$ and $q=5$. For each of the cases $n=20$ was assumed. Accordingly, the curve created by the scriber mounted at the 20$^{th}$ wheel for $q=2.5$ is shown in Fig. 2.

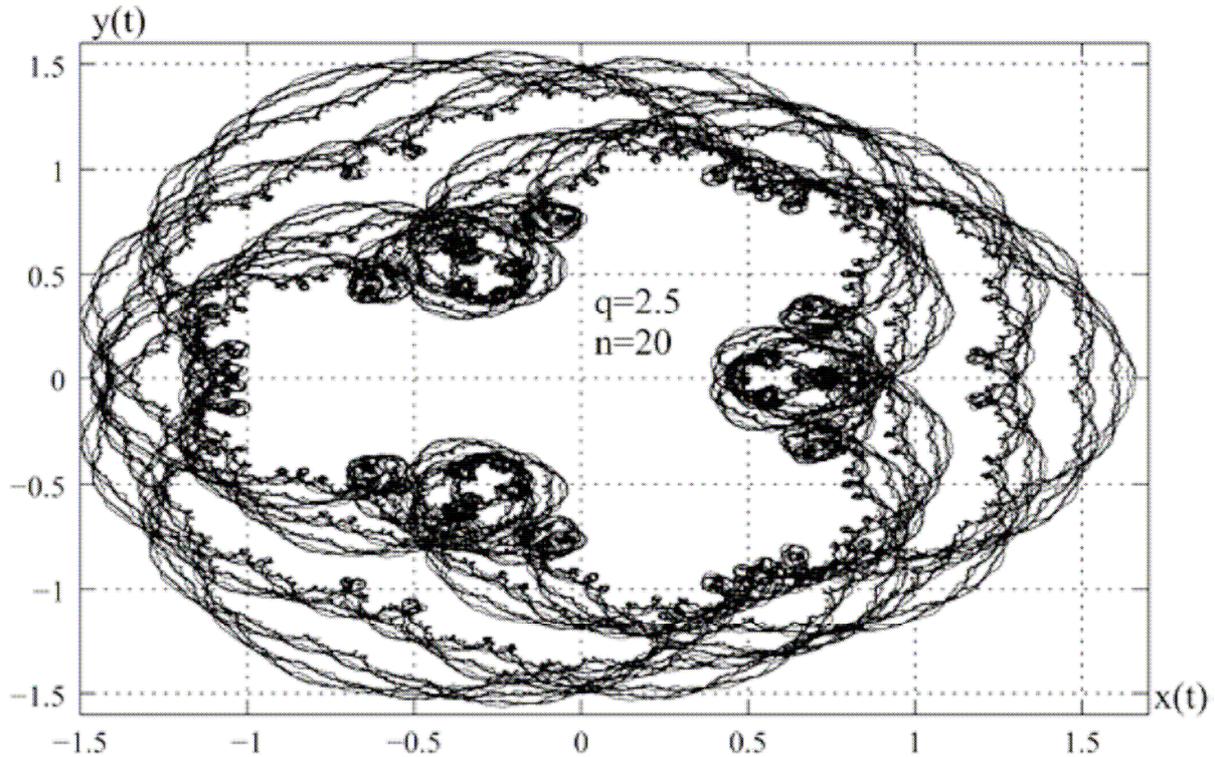

**Fig. 2.** Phase diagram of the scriber on the phase plane. $q = 2.5$, $n = 20$.

The presented geometrical structure is very complex, but the images of its fragments point to a fractal form (Fig. 3 and 4) [1], [3-4].

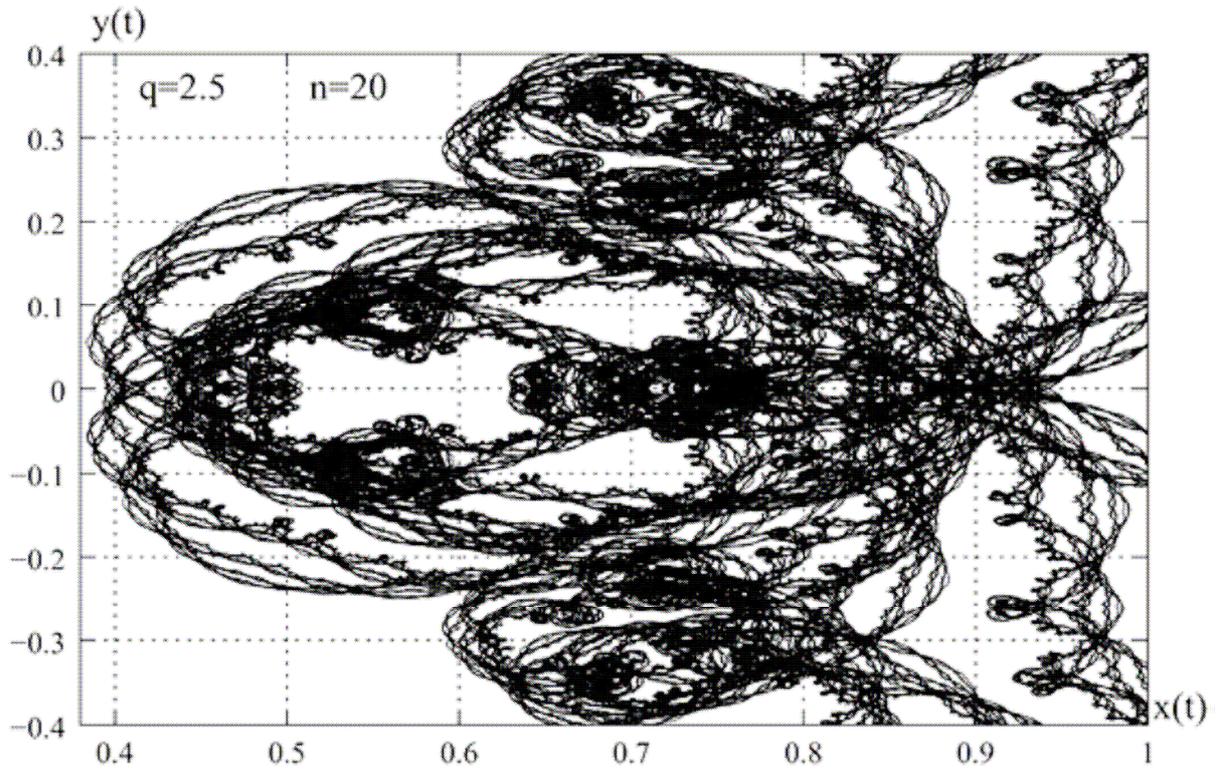

**Fig. 3.** Fragment of Figure 2.

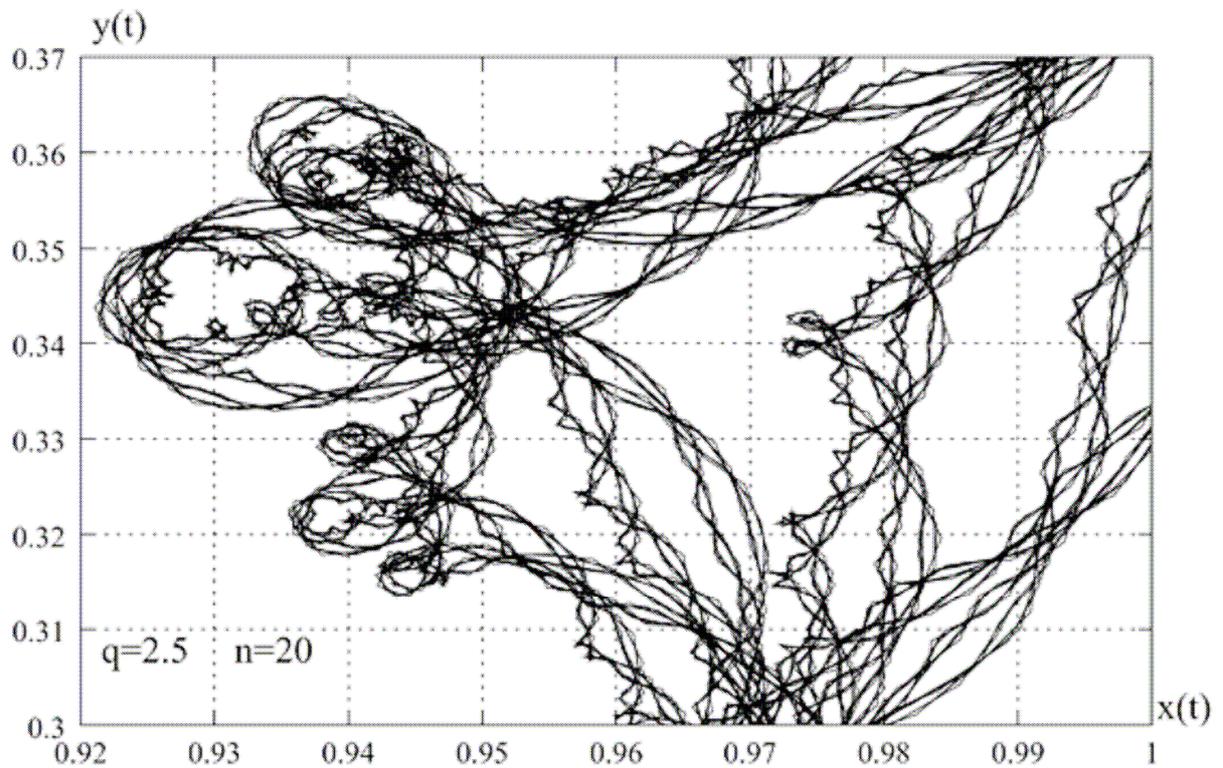

**Fig. 4.** Fragment of Figure 3.

Coordinates $x$ and $y$ are normalized in accordance with the dependence: $x = x'/r_0$, $y = y'/r_0$. Periodic changes of the length of radius-vector $|R|$ are illustrated in Fig. 5, whereas in Fig. 6 the speed of $\left|\dfrac{dR}{dt}\right|$ changes is shown.

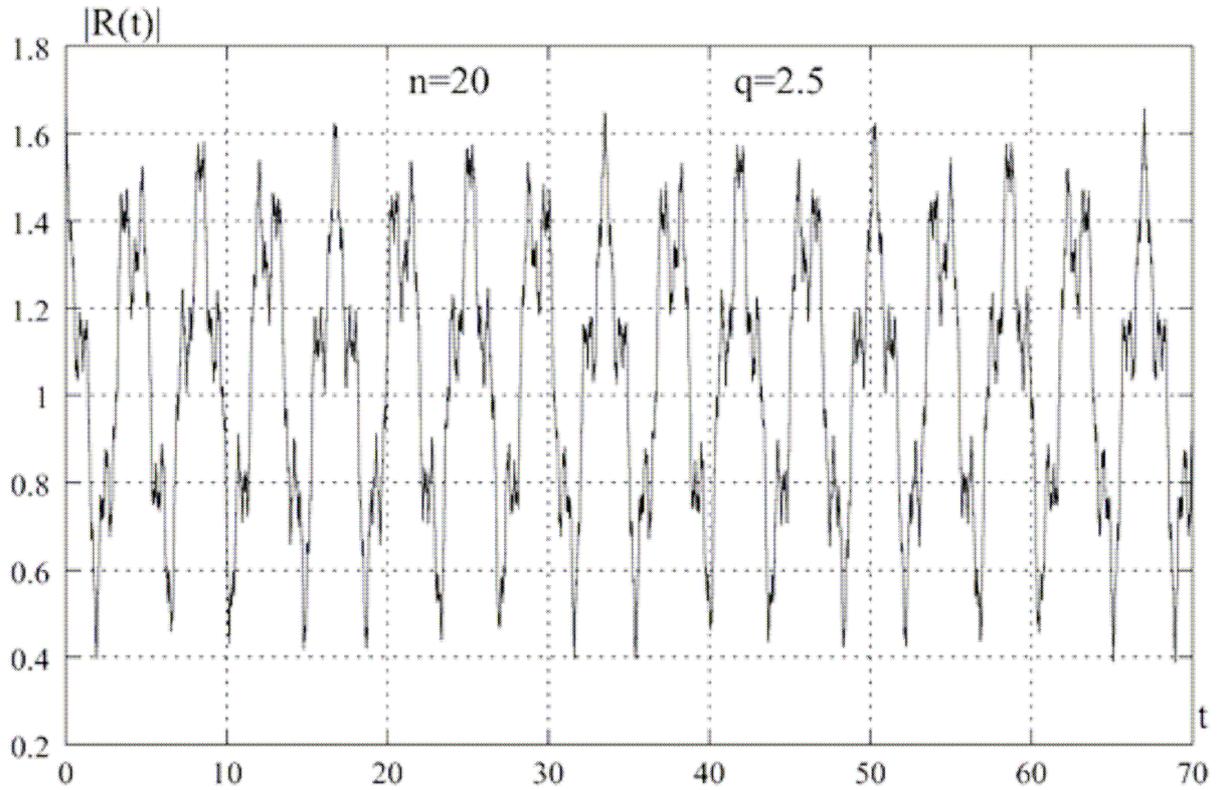

**Fig. 5.** The dependence between the radius-vector and the angle of rotation of the biggest wheel. *q* = 2.5, *n* = 20.

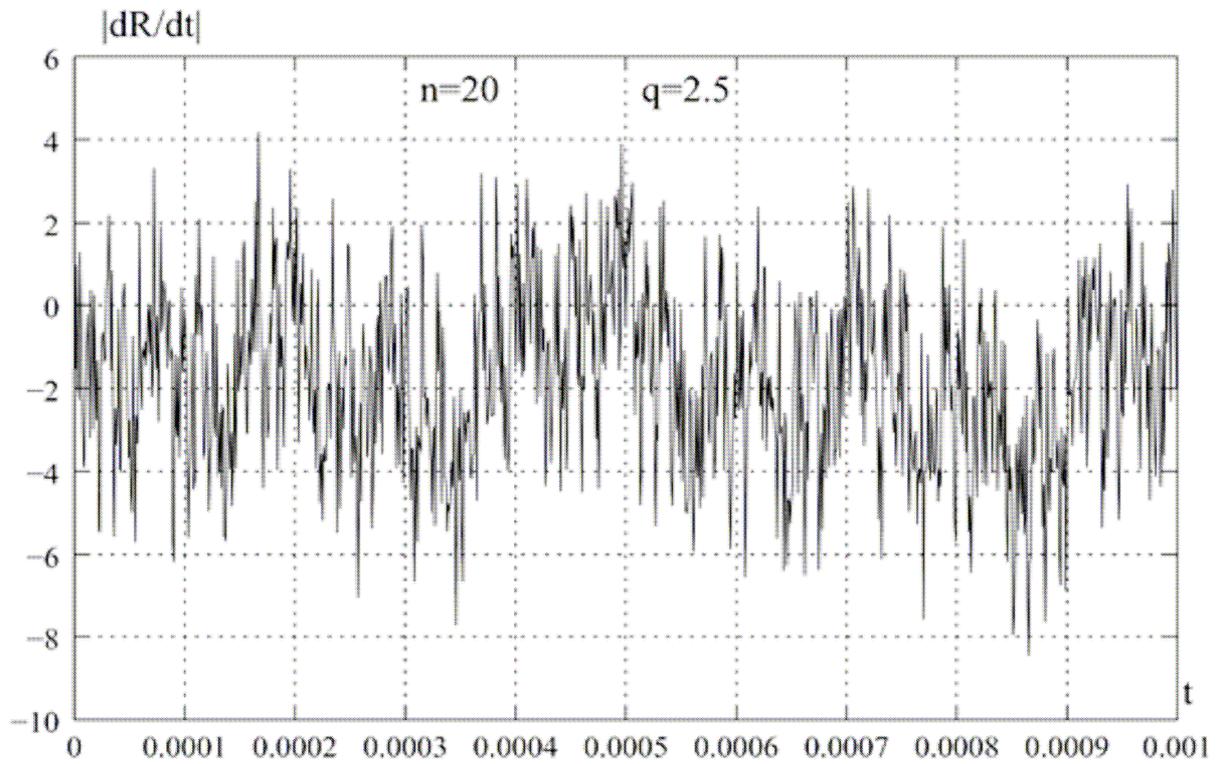

**Fig. 6.** Dependence between the rate of the changes of the radius-vector and the angle of rotation of the biggest wheel. *q* = 2.5, *n* = 20.

The graphs certify the continuity of *R(t)* and the discontinuity of its derivative in practice. In an extreme case, for each value of angle *t* a rapid change of speed occurs. The zigzags visible on the phase plane are a consequence of the above-mentioned discontinuity (Fig. 7).

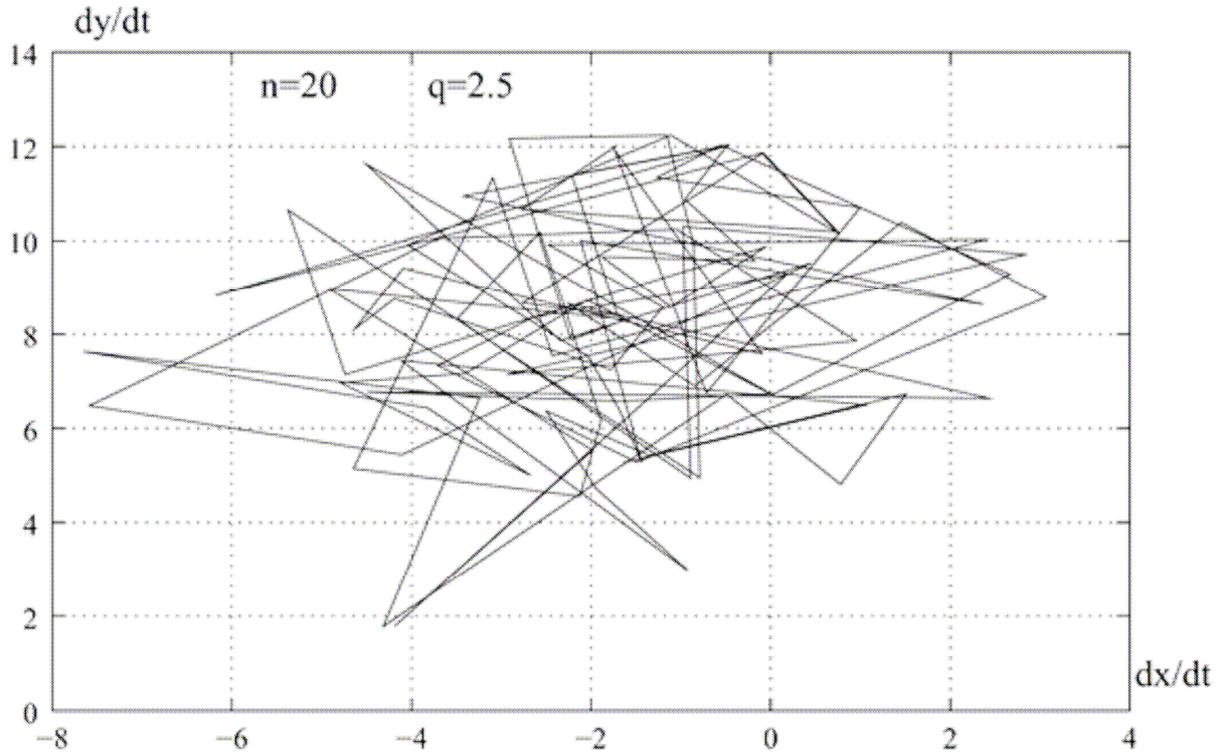

**Fig. 7.** Phase trajectory of the rate of changes in the coordinates of the scriber.
$q = 2.5, n = 20$.

When $n = \infty$ the zigzags occur at each point of the plan. Such trajectory resembles Brownian motion, where, as commonly known, the motion of basic particles is also zigzag-like at each point of the surface [7-8].

It should be emphasised that even though the motion of the scriber is not of a chaotic nature, yet it is very sensitive to changes in ratio $\frac{\varphi_{k+1}}{\varphi_k}$. Thus, in Fig.8 the thick line represents the case when $q = 2.5$, whereas the thin line the case of $q_\varphi = \frac{\varphi_{k+1}}{\varphi_k} = 2.499$.

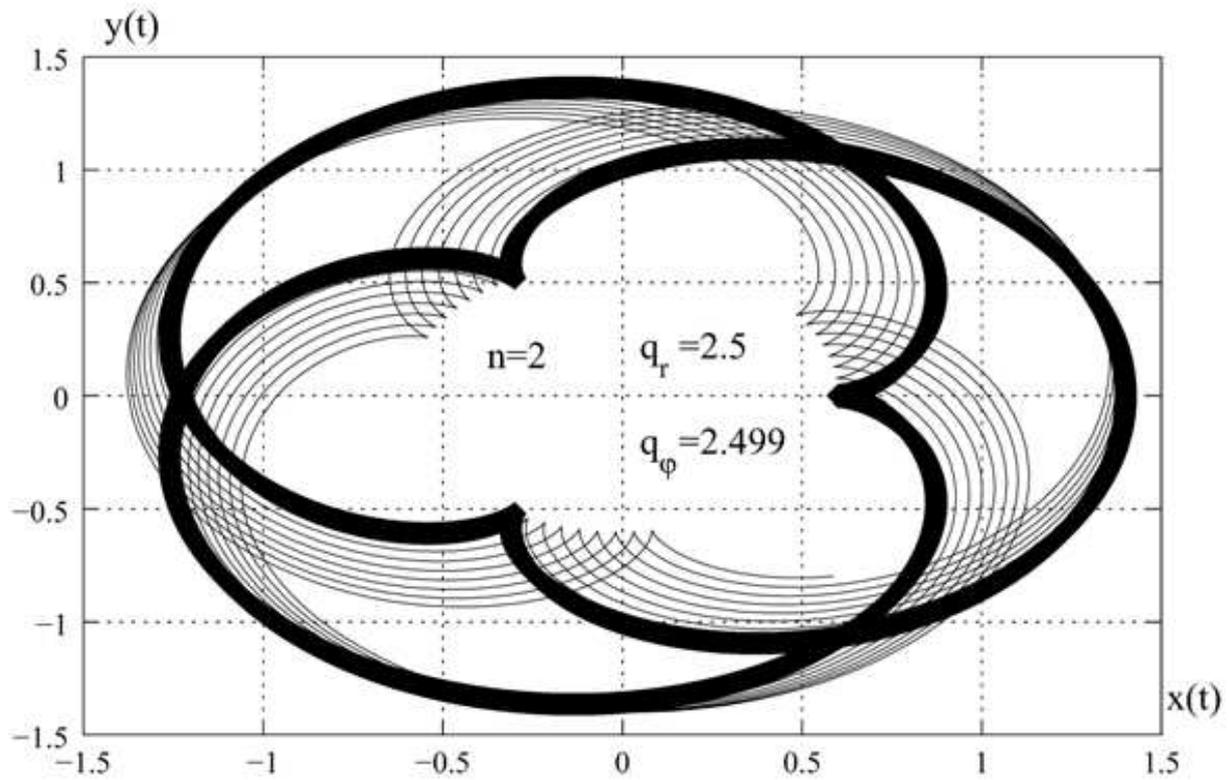

**Fig. 8.** The sensitivity of the system to changes in the angular speed ratios. The thick line: $q = 2.5$. The thin line: $qr = 2.5$, $q\phi = 2.499$.

Both graphs refer to $n=2$. Yet, the change of the value of $q_r$, given the same $q_\varphi$ value, does not affect significant changes in the course of the trajectory.

In Fig.9 the phase trajectory of the scriber motion was shown for $q=5$.

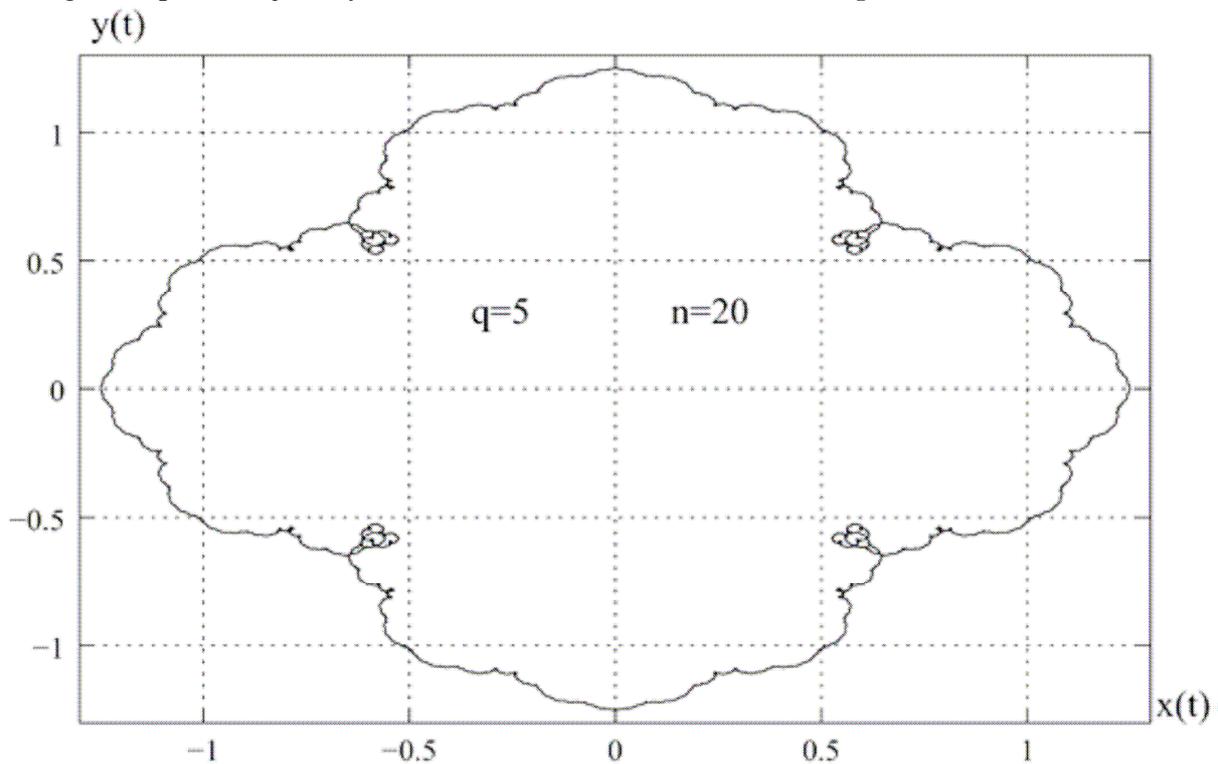

**Fig. 9.** Phase diagram of the motion of the scriber on the plane. $q = 5, n = 20$.

It may be observed that this motion is not as complex as in the previously discussed case, nonetheless, the derived graph has a fractal nature, as substantiated by the fragments presented in Figures 10 and 11. After scaling, the Figures are identical.

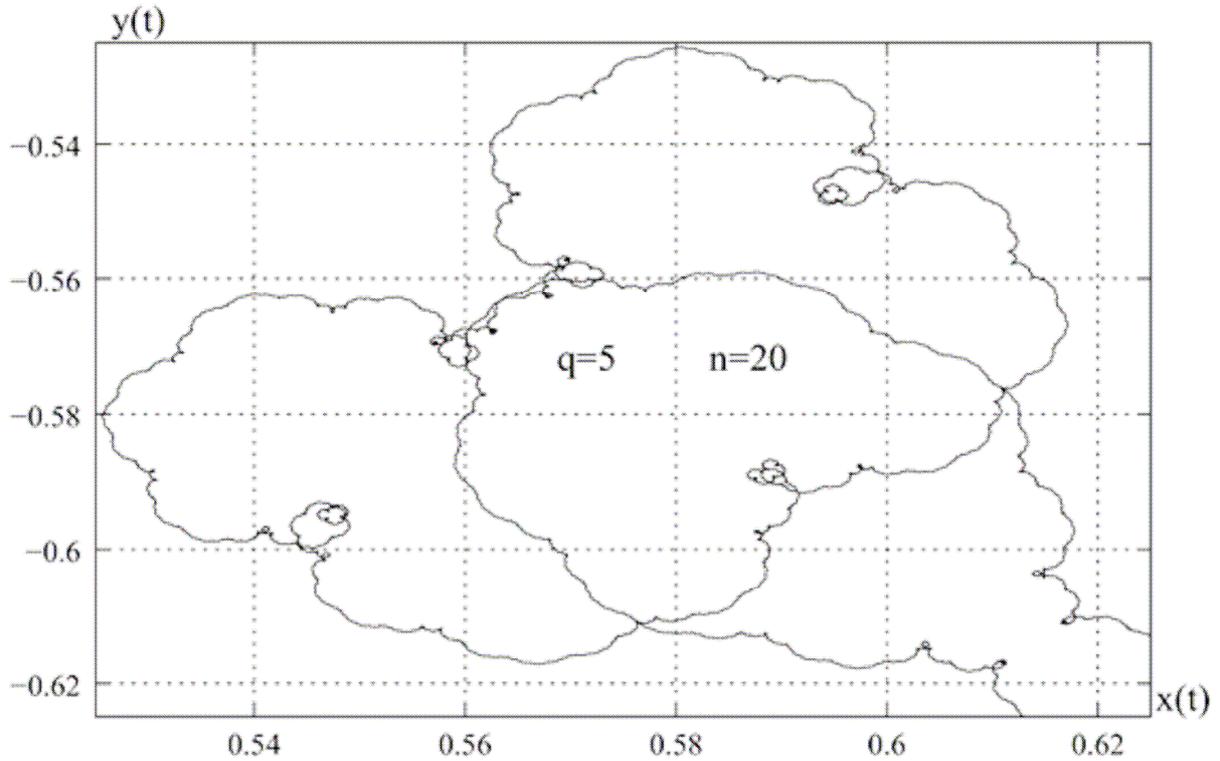

**Fig. 10.** Fragment of Figure 9.

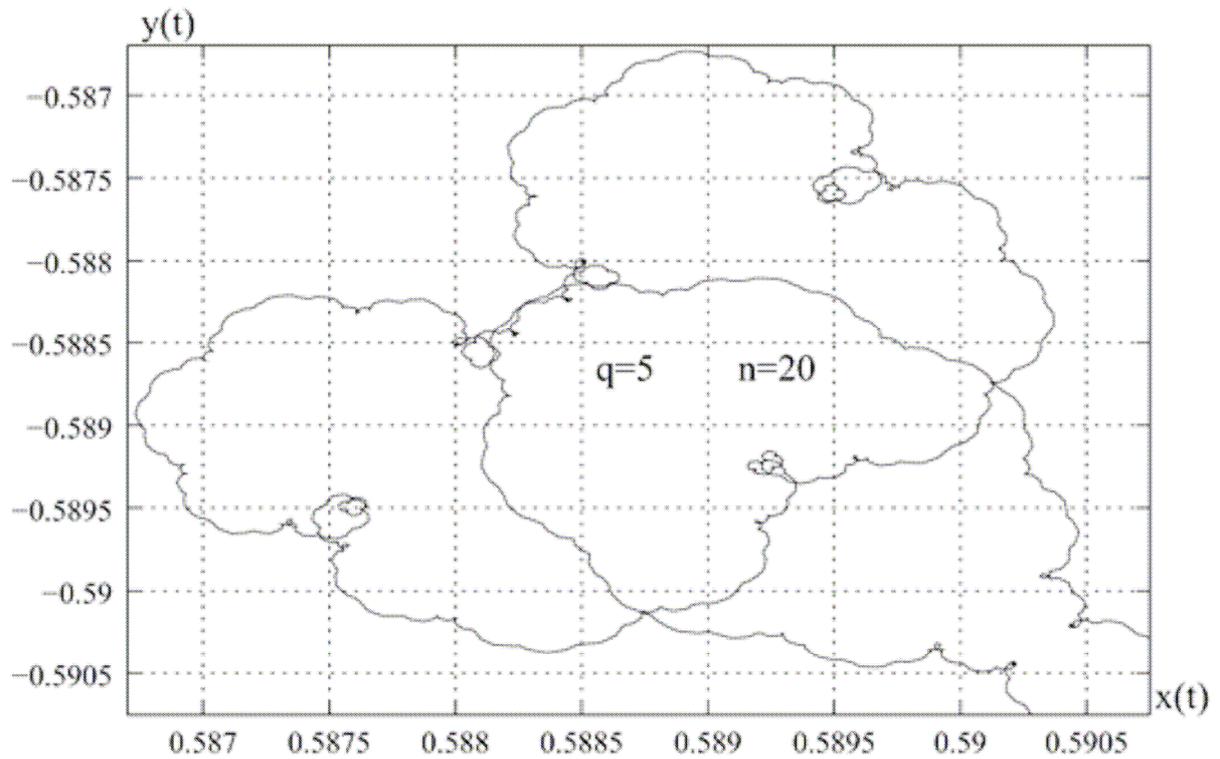

**Fig. 11.** Fragment of Figure 10.

In Fig.12 changes in the information system entropy are indicated, depending on $q$. They also have a fractal nature. The values of the entropies were calculated from the following formula:

$$E = -\sum_{i=1}^{N} p_i \lg_2 p_i \qquad (9)$$

where $p_i$ is the probability of the occurrence of a definite value of the length of radius $|R|$ [6]. As seen in the above graph, most information is rendered by the system for $q=1.8$.

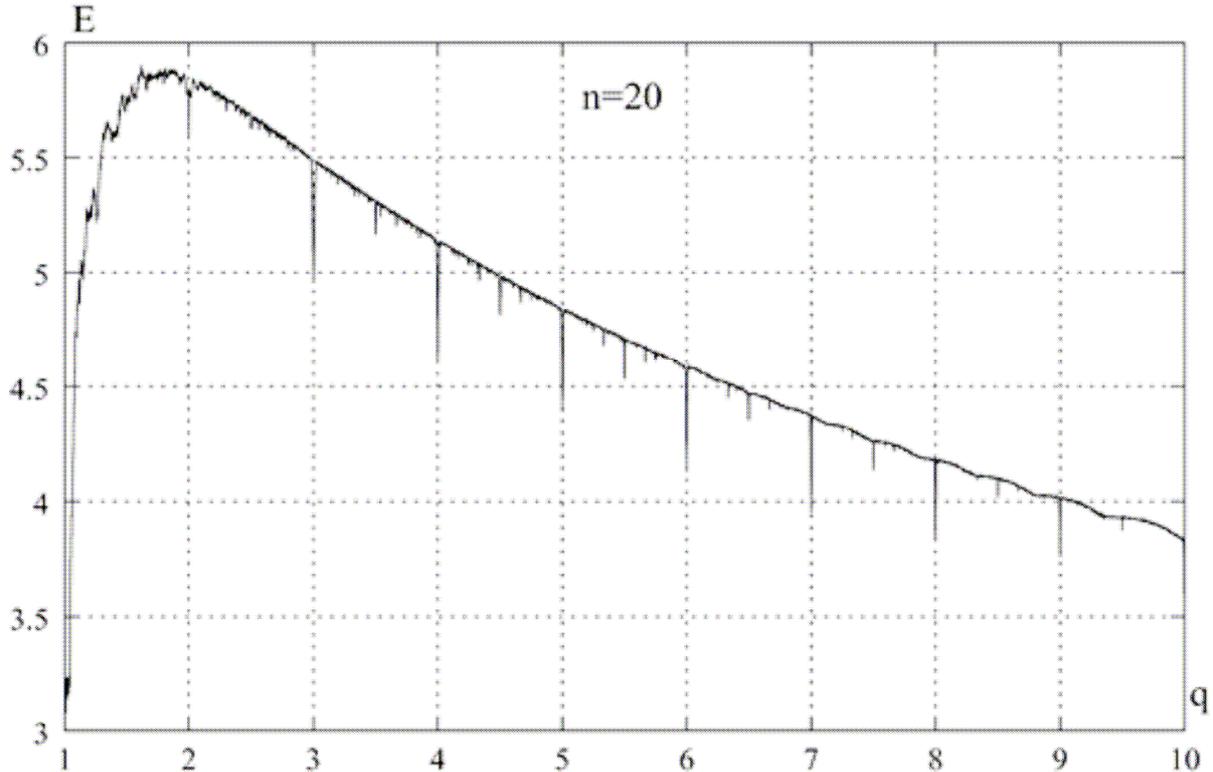

**Fig. 12.** Dependence between the information system entropy and $q$ ratio.

## 3. Concluding remarks

A mathematical and numerical analysis conducted within the framework of the paper was focused on the dynamics of a mechanical system consisting of $n$ wheels rotating at constant speed around their axes. The smaller wheels were mounted on the rim of bigger wheels- see Figure 1. It was assumed that the ratios between the radii and the angles of rotation of the successive wheels are constant. There is a scriber mounted on the rim of the smallest wheel. In the outcome of the analysis, very complex fractal graphs were derived. The behavior of the speed of the scriber is very interesting, as its components, represented on the phase plane, resemble Brownian motion (Fig.7).

Although the changes in the length of the radius vector have a continuous character (Fig. 5), yet the course of the changes is not even. With an infinite number of the wheels the roughness of the graph occurs at each of its points. A spatial form of this phenomenon is shown in Fig.13.

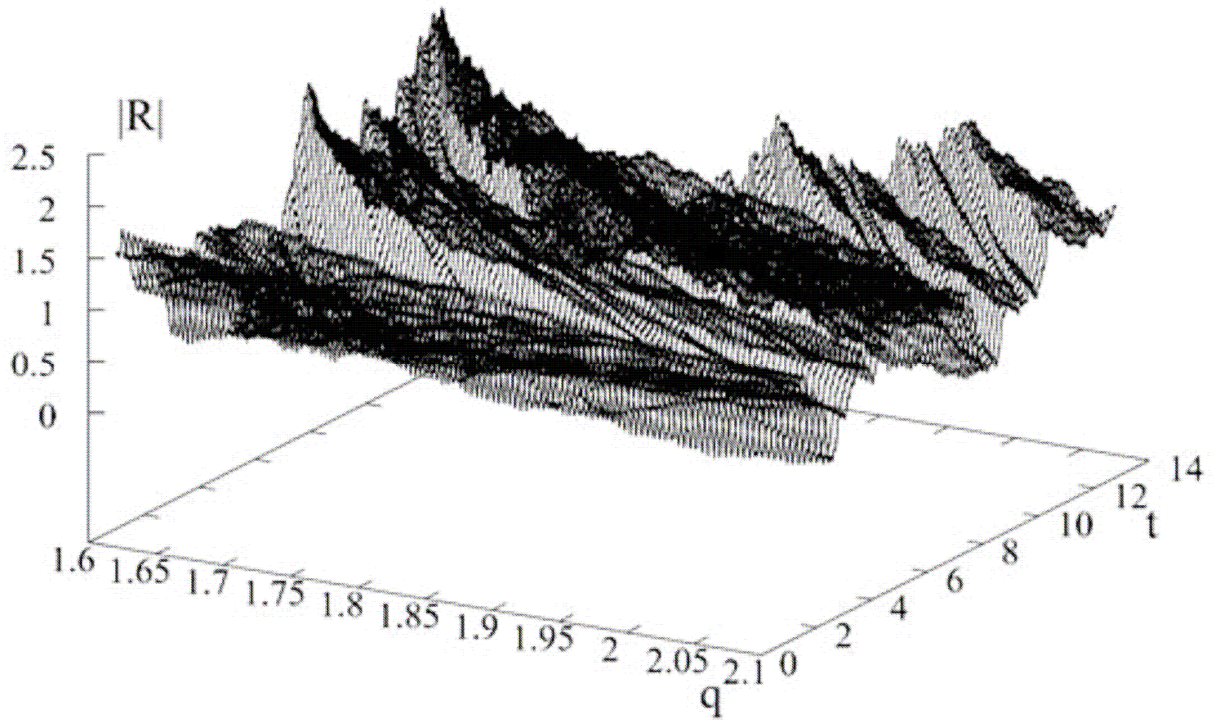

**Fig. 13.** Visualization space of the changes of the radius-vector.

It is difficult to refrain from the impression that natural biological structures have the same character. It is sufficient to have a look at afforested mountain slopes, rocky or non-stony. Therefore, a conclusion may be drawn that in the natural environment we can also witness the phenomenon of common roughness.

## Notations

| | |
|---|---|
| $n$ | number of wheels |
| $N$ | observations horizon |
| $p$ | probability |
| $r$ | wheel radius |
| $R$ | radius-vector |
| $t$ | angle of rotation (angular displacement) of the biggest wheel |
| $x,y$ | scriber position coordinates |
| $\varphi$ | angle |
| $q$ | ratio of radii and angular speeds |

*Subscripts*

| | |
|---|---|
| $k$ | wheel number |
| $r$ | refers to the wheel radius |
| $\varphi$ | refers to the angle of rotation of the wheel |
| $0$ | refers to the biggest wheel |

**References**
[1] Berezowski M. Fractal solutions of recirculation tubular chemical reactors. Chaos, Solitons & Fractals, 2003; 16, 1-12.
[2] Berezowski M. Analysis of dynamic solutions of complex delay differential equation exemplified by the extended Mandelbrot's equation. Chaos, Solitons & Fractals, 2005; 24: 1399 – 1404.
[3] Berezowski M. Fractal character of basin boundaries in a tubular chemical reactor with mass recycle. Chemical Engineering Science, 2006; 61: 1342-1345.
[4] Berezowski M. & Bizon K. Fractal structure of iterative time profiles of a tubular chemical reactor with recycle. Chaos, Solitons & Fractals, 2006; 28: 1046-1054.
[5] Berry M. V. & Lewis, Z. V. On the Weierstrass-Mandelbrot Function. Proc. Roy. Soc. London Ser., 1980; A 370: 459-484.
[6] Du Bois-Reymond P. Versuch einer Classification der willk urlichen Functionen reeller Argumente nach ihren Aenderungen in den kleinsten Intervallen. J. Reine Angew. Math., 1875; 79: 21–37.
[7] Einstein A. Über die von der molekularkinetischen Theorie der Wärme geforderte Bewegung von in ruhenden Flüssigkeiten suspendierten Teilchen. Annalen der Physik, 1905; 17: 549–560.
[8] Smoluchowski M. Zur kinetischen Theorie der Brownschen Molekularbewegung und der Suspensionen. Annalen der Physik, 1906; 21: 756–780.
[9] Song L., Xu S. & Yang J. Dynamical models of happiness with fractional order. Communications in Nonlinear Science & Numerical Simulation, 2010; 15: 616 – 628.
[10] Tang Y. & Fang J. Synchronization of *N*-coupled fractional-order chaotic systems with ring connection. Communications in Nonlinear Science & Numerical Simulation, 2010; 15: 401 - 412.